\documentclass[aapm,mph,reprint,graphicx]{revtex4-1}
\usepackage[pdftex]{graphicx}
\usepackage[colorlinks=true,linkcolor=blue,citecolor=blue]{hyperref}

\DeclareGraphicsExtensions{.pdf,.jpeg,.png}
\usepackage{amsmath,amssymb,tikz}

\newcommand{\xvec}{\mathbf{x}}
\newcommand{\bvec}{\mathbf{b}}
\newcommand{\muvec}{\boldsymbol \mu}
\newcommand{\tvec}{\mathbf{t}}
\newcommand{\pvec}{\mathbf{p}}
\newcommand{\cvec}{\mathbf{c}}

\usepackage{url}

\begin{document}

\def\papertitle{Technical Note: Convergence analysis of a polyenergetic SART algorithm}
\title{\papertitle}

\author{T.~Humphries}
\email[Email: ]{humphrit@science.oregonstate.edu}
\affiliation{Department of Mathematics, Oregon State University}

\date{\today}

\begin{abstract}	
\textbf{Purpose:} We analyze a recently proposed polyenergetic version of the simultaneous algebraic reconstruction technique (SART). This algorithm, denoted pSART, replaces the monoenergetic forward projection operation used by SART with a post-log, polyenergetic forward projection, while leaving the rest of the algorithm unchanged. While the proposed algorithm provides good results empirically, convergence of the algorithm was not established mathematically  in the original paper. 

\textbf{Methods:} We analyze pSART as a nonlinear fixed point iteration by explicitly computing the Jacobian of the iteration. A necessary condition for convergence is that the spectral radius of the Jacobian, evaluated at the fixed point, is less than one. A short proof of convergence for SART is also provided as a basis for comparison.

\textbf{Results:} We show that the pSART algorithm is not guaranteed to converge, in general. The Jacobian of the iteration depends on several factors, including the system matrix and how one models the energy dependence of the linear attenuation coefficient. We provide a simple numerical example that shows that the spectral radius of the Jacobian matrix is not guaranteed to be less than one. A second set of numerical experiments using realistic CT system matrices, however, indicates that conditions for convergence are likely to be satisfied in practice.

\textbf{Conclusion:}  Although pSART is not mathematically guaranteed to converge, our numerical experiments indicate that it will tend to converge at roughly the same rate as SART for system matrices of the type encountered in CT imaging. Thus we conclude that the algorithm is still a useful method for reconstruction of polyenergetic CT data.
\end{abstract}

\maketitle


\section{Introduction}
The algebraic reconstruction technique (ART) or Kaczmarz' method~\cite{K37, GBH70} is a well-known method for approximately solving systems of linear equations,
\begin{equation}\label{E:Axb}
A\xvec = \bvec,
\end{equation}
where $\xvec$ is a column vector of size $n$, $\bvec$ is a column vector of size $m$, and $A$ is $m \times n$. The method has a long history in CT image reconstruction, in which $\bvec$ represents the post-log measured projection data, $\xvec$ represents the image to be reconstructed, and the $(i,j)$th entry of $A$ represents the length or area of intersection of the $i$th ray with the $j$th pixel of the image. The values in $\xvec$ represent the averaged linear attenuation coefficient (LAC) within each pixel, and are usually constrained to be positive. Since the model is linear, it is implicitly assumed that the X-ray beam used to generate the data is monoenergetic, and the reconstructed values correspond to the LAC of the tissue at that energy.

Beginning from an initial estimate $\xvec^{(0)}$, ART iteratively projects the current iterate, $\xvec^{(k)}$, onto the hyperplane defined by one of the $m$ equations defined by (\ref{E:Axb}). The ART iteration is given by
\begin{align}
 \xvec^{(k+1)} &= \xvec^{(k)} - \frac{\langle a_i,  \xvec^{(k)} \rangle - b_i}{\langle a_i , a_i \rangle} a_i^T, \label{E:ART}
\end{align}
where $a_i$ is the $i$th row of $A$, $\langle \cdot, \cdot \rangle$ denotes the dot product, and $i$ is chosen to be $(k \mod m) + 1$ in the classical version of the algorithm. Some references consider this step to be a sub-iteration, and a full iteration of ART to consist of successively applying this iteration for all $m$ equations. If the system is consistent, the iteration is guaranteed to converge to a solution, while in the inconsistent case the sub-iterations converge to a limit cycle~\cite{T71}.

The simultaneous ART (SART) method~\cite{AK84,KS01} is a variant of ART in which the corrections generated by the ART sub-iterations (\ref{E:ART}) are combined and applied simultaneously. The iteration can be expressed concisely in terms of matrix operations as 
\begin{equation}
 \xvec^{(k+1)} = \xvec^{(k)} - D A^T M\left( A  \xvec^{(k)} - \bvec \right), \label{E:SART}
\end{equation}
where $D$ and $M$ are diagonal matrices, with
\begin{align}
D_{jj} &= \frac{1}{\beta_j}, ~~\beta_j = \sum_{k=1}^m |a_{kj}|, j = 1 \dots n \notag\\
M_{ii} &= \frac{1}{\gamma_i}, ~~\gamma_i = \sum_{k=1}^n |a_{ik}|, i = 1 \dots m. \label{E:diagmatrices} 
\end{align}
In other words, $\beta_j$ is the 1-norm of the $j$th column of $A$, and $\gamma_i$ is the 1-norm of the $i$th row. 

The primary advantage of SART over ART is that SART is less sensitive to noisy data~\cite{AK84}. It should be noted that the iteration~(\ref{E:SART}) computes an update using projection data corresponding to all views simultaneously, while the original SART algorithm~\cite{AK84} computes a sequence of updates using projection data corresponding to only one view of the object at one time, to accelerate convergence. For simplicity of analysis we will only consider equation~(\ref{E:SART}), as this was the used as the basis for the polyenergetic algorithm that we will analyze. Like ART, SART has been proven to converge~\cite{CE02,JW03}; see also the Appendix to this paper. 

As mentioned, both ART and SART implicitly assume that the data are generated from a monoenergetic X-ray beam, i.e. that 
\begin{equation}
p_i = I \exp \left( - \langle a_i, \xvec \rangle \right), \label{E:monoCT1}
\end{equation}
where $p_i$ is the measured intensity of the $i$th beam, and $I$ is the blank scan intensity (assumed to be independent of $i$). Taking the log of the data and rearranging terms then gives a linear system equivalent to~(\ref{E:Axb}).
\begin{equation}
- \ln \frac{p_i}{I} = \langle a_i, \xvec \rangle. \label{E:monoCT2}
\end{equation}
In practice, however, the X-rays generated by clinical CT hardware are usually polyenergetic.  A typical model for polyenergetic X-ray measurements is
\begin{equation}
p_i = \int I (\varepsilon) \exp \left( - \langle  a_i, \muvec (\varepsilon) \rangle \right) \: d \varepsilon, \label{E:polyCT1}
\end{equation}
where $\varepsilon$ refers to the energy of an incident x-ray, $I(\varepsilon)$ is the initial intensity of the beam corresponding to that energy (i.e. the beam's spectrum) and $\muvec(\varepsilon)$ is a vector of attenuation coefficients, whose values depend on X-ray energy. This system of equations can no longer be linearized, and it is well-known that reconstructing an image from this data using conventional means (such as ART, SART, or filtered back projection) produces an image containing beam hardening artifacts~\cite{BD76}. This motivates the need for polyenergetic iterative reconstruction algorithms, e.g. \cite{HMDJ00,YWBYN00,DNDMS01, EF02, EF03,VVDB11,Rezvani12,HF14,LS14a, LS14b}.

The polyenergetic SART (pSART) algorithm~\cite{LS14b} was recently proposed as one such reconstruction technique. In this approach, the beam spectrum is discretized into $h$ energy levels $\varepsilon_h$, with weighting terms $I_h$ computed to approximate the continuous spectrum. The vector-valued function $\muvec(\varepsilon)$ is then modeled as a function of a vector $\tvec$, representing the attenuation map of the object at a reference energy, $\varepsilon_0$, which was chosen to be 70 keV. This function makes use of tabulated energy-dependent LAC values for some suitable reference materials, such as air, fat, breast, soft tissue and bone. Letting $t_j$ denote the LAC value for pixel $j$ at the reference energy, the LAC of that pixel for all other energies $\varepsilon$ is then given by
\begin{equation}
\mu(t_j, \varepsilon) = \frac{ [\mu_{k+1}(\varepsilon_0) - t_j] \mu_{k}(\varepsilon) +  [t_j - \mu_{k}(\varepsilon_0) ] \mu_{k+1}(\varepsilon)}{\mu_{k+1}(\varepsilon_0) - \mu_{k}(\varepsilon_0)}, \label{E:interp}
\end{equation}
where $\mu_k(\varepsilon)$ and $\mu_{k+1}(\varepsilon)$ are the tabulated, energy-dependent LAC functions for the two base materials with LAC values adjacent to $t_j$ at the reference energy. So for instance, if the value of $t_j$ is between the LAC for soft tissue and the LAC for bone at the reference energy, then its LAC at all other energies is obtained by linear interpolation between the corresponding values for bone and soft tissue, with the weighting determined by the values at the reference energy. 

One can then define a polyenergetic forward projection operator, $\mathcal{P}: \mathbb{R}^n \to \mathbb{R}^m$, which acts on $\tvec$:
\begin{equation}
[\mathcal{P} (\tvec)]_i = \sum_h I_h \exp \left( - \langle a_i, \muvec (\tvec, \varepsilon_h) \rangle \right). \label{E:project_poly}
\end{equation}
The pSART iteration is then defined as
\begin{equation}
 \tvec^{(k+1)} = \tvec^{(k)} - D A^T M\left( -\ln  \left[\mathcal{P}  \left( \tvec^{(k)} \right) \right] +\ln (\pvec) \right), \label{E:pSART}
\end{equation}
with $D$ and $M$ defined as in~(\ref{E:diagmatrices}). Equations~(\ref{E:project_poly}) and (\ref{E:pSART}) are equivalent to equations (14) and (15) in Ref.~\onlinecite{LS14b}, although our notation is somewhat different. The only difference from the SART iteration (\ref{E:SART}) is that the log of the monoenergetic forward projection has been replaced by the log of the polyenergetic forward projection. The algorithm produces a single attenuation map of the object, $\tvec$, with LAC values corresponding to the reference energy. 

In Ref.~\onlinecite{LS14b} the authors state that this modification solves the problem of inconsistency between the polyenergetic data and the monoenergetic model implicitly assumed by the conventional SART approach. While it is clear that a vector $\tvec$ satisfying $\pvec = \mathcal{P} (\tvec)$ is a fixed point of this iteration, this does not guarantee convergence of the algorithm. Experimental results indicated that the method is effective, however. In the next section we analyze the convergence of pSART as a fixed point iteration.

\section{Convergence of pSART}

We first consider the SART iteration~(\ref{E:SART}). The iteration can be written in the form
\begin{equation}
\xvec^{(k+1)} = T\xvec^{(k)} + \cvec, \label{E:SART2}
\end{equation}
where $T = I - DA^TMA$ and $\cvec = DA^TM\bvec$. One can show that the spectral radius of $T$ (the magnitude of its largest eigenvalue), denoted $\rho(T)$, is strictly less than 1  (see Appendix). This guarantees convergence of the algorithm to a solution, if the system is consistent. It has been proven that if no exact solution exists, SART converges to a weighted least-squares solution~\cite{CE02, JW03}.

The pSART algorithm is a nonlinear fixed point iteration. We write the iteration as
\begin{align}
\tvec^{(k+1)} &=  \tvec^{(k)} - DA^TM f\left(\tvec^{(k)}\right), \label{E:pSART2} \\
		   &\equiv F\left(\tvec^{(k)}\right) 
\end{align}
where $f(\tvec) = -\ln  \left[\mathcal{P}  (\tvec) \right] +\ln (\pvec) $. Assuming that a solution to the system of nonlinear equations exists (which corresponds to a fixed point of the iteration), a necessary condition for convergence is that the spectral radius of the Jacobian matrix of $F$, $\rho(J_F)$, must be less than one when evaluated at the solution. Note that this does not guarantee that the algorithm converges to the solution from any starting point; simply that it cannot converge if the condition does not hold. Let $\tvec^*$ denote a solution to the nonlinear system of equations $\pvec_i = [ \mathcal{P} (\tvec)]_i$. A straightforward calculation gives
\begin{equation}
J_F\left(\tvec^*\right) = I - DA^TM J_f\left(\tvec^*\right),\label{E:jac1}
\end{equation}
where the $(i,j)$th element of the $m \times n$ Jacobian matrix $J_f(\tvec)$ is given by
\begin{equation}
\frac{ \partial f_i}{\partial t_j}(\tvec) = \frac{a_{ij}}{[\mathcal{P} (\tvec)]_i } \left[ \sum_h I_h \exp \left( - \langle  a_i, \muvec (\tvec, \varepsilon_h) \rangle \right) \frac{ \partial \mu}{\partial t} (t_j, \varepsilon_h) \right].\label{E:jac2}
\end{equation}
It follows from~(\ref{E:interp}) that
\begin{equation}
\frac{ \partial \mu}{\partial t} (t, \varepsilon) = \frac{ \mu_{k+1}(\varepsilon) - \mu_k (\varepsilon) }{ \mu_{k+1}(\varepsilon_0) - \mu_k (\varepsilon_0) },
\end{equation}
where $\mu_{k+1}$ and $\mu_k$ again refer to the base material LACs that are adjacent to $t$ at the reference energy $\varepsilon_0$. This implies that for a fixed value of $\varepsilon$, this partial derivative is piecewise constant, with discontinuities when $t$ is equal to one of the reference LACs at $\varepsilon_0$.

A general analysis of the spectral radius of $J_F(\tvec^*)$ is difficult as it has a complicated dependence on the system matrix $A$, the spectrum $I_h$, and the choice of base materials. We now show with a numerical experiment that one cannot guarantee that $\rho(J_F(\tvec^*)) < 1$, in general. 

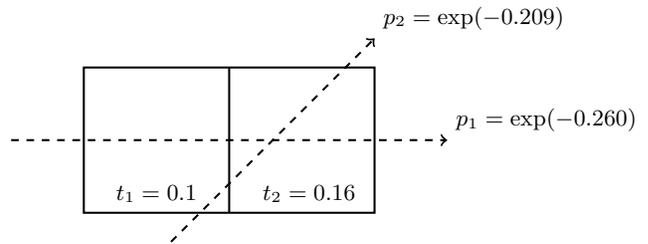
\begin{figure}
\begin{center}
\resizebox{\linewidth}{!}{
\begin{tikzpicture}

\draw[black,thick] (0,0) rectangle (4,2);
\draw[black,thick] (2,0) -- (2,2);
\draw[black,thick,dashed,->] (-1,1) -- (5,1) node[anchor=south west]{$p_1 = \exp(-0.260)$};
\draw[black,thick,dashed,->] (1.2,-0.4) -- (4,  2.4) node[anchor=south west]{$p_2 = \exp(-0.209)$};
\node[anchor = north] (n1) at (1,0.5) {$t_1 = 0.1$};
\node[anchor = north] (n2) at (3.1,0.5) {$t_2 = 0.16$};

\end{tikzpicture}
}
\end{center}
\caption{Example problem consisting of two pixels.}\label{F:exp1}
\end{figure}

\begin{figure}
\includegraphics*[width=0.45\linewidth]{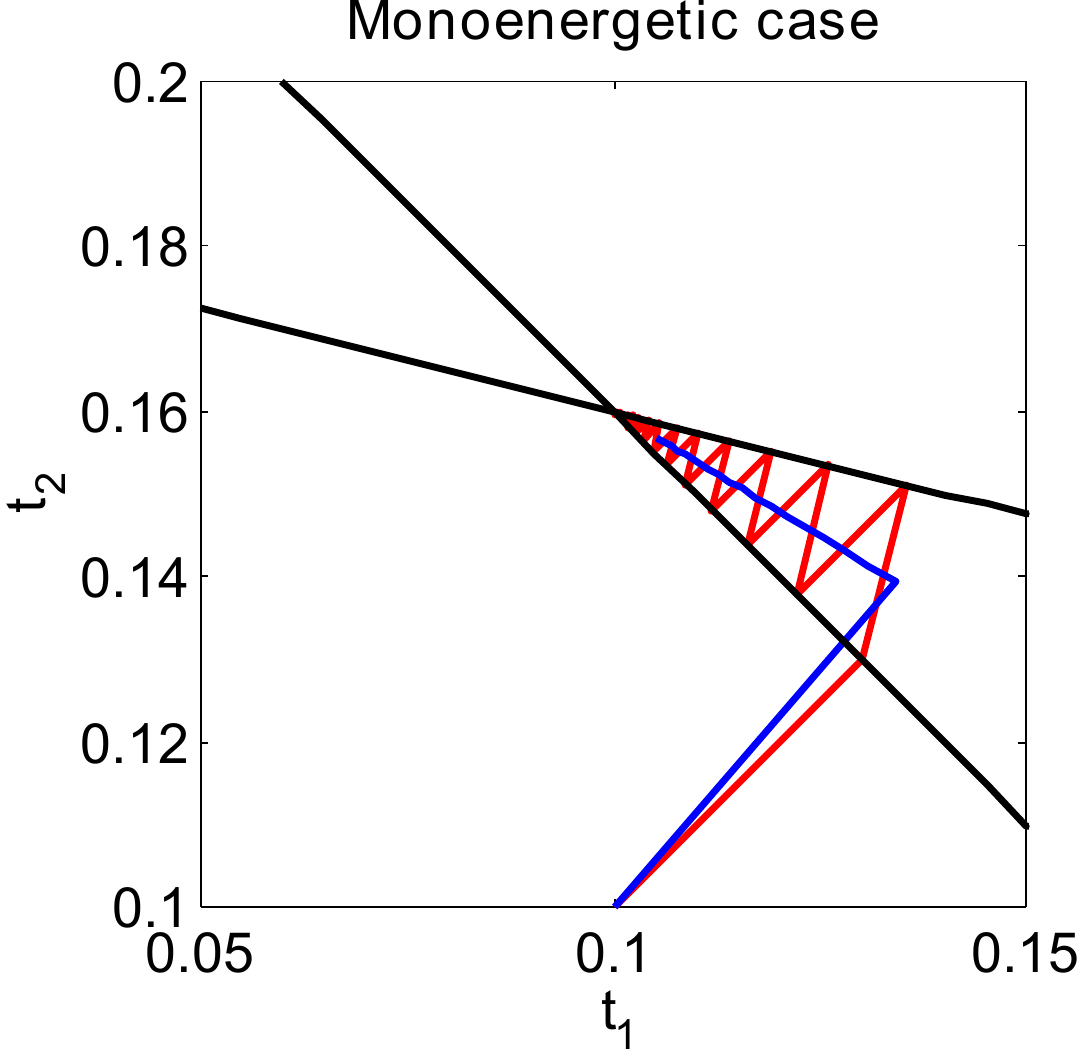}
\caption{Convergence of ART \& SART for the monoenergetic experiment. The black lines correspond to the system of equations (\ref{E:eqs1}), red line correspond to the sequence of iterates generated by ART, blue line to the sequence of iterates generated by SART.}\label{F:ART_mono}
\end{figure}

We consider a simple example where $m = n = 2$, illustrated in Fig.~\ref{F:exp1}. The object consists of two pixels of size 1$\times$1 cm, with LACs of $t_1 = 0.1$ cm$^{-1}$ and $t_2 = 0.16$ cm$^{-1}$ at a reference energy of 70 keV. We first consider the case of a monoenergetic 70 keV beam with intensity $I = 1$. The first beam travels through both pixels horizontally, while the second beam has a length of intersection of roughly 0.28 cm with the first pixel and 1.13 cm with the second pixel. After taking logarithms, we obtain the $2\times 2$ linear system of equations
\begin{equation}
\begin{bmatrix}
1 & 1 \\ 0.28 &1.13
\end{bmatrix}
\begin{bmatrix} t_1 \\ t_2 \end{bmatrix}
= \begin{bmatrix}0.260 \\ 0.209 \end{bmatrix}, \label{E:eqs1}
\end{equation}
which can be solved using ART or SART. Fig.~\ref{F:ART_mono} shows the progression of both ART and SART for this system, starting from an initial guess of zero. One can see that ART sequentially projects onto the two lines defining the equations, while the iterates generated by SART follow a path between the two equations.

\begin{figure}
\begin{tabular}{cc}
\includegraphics[width=0.45 \linewidth]{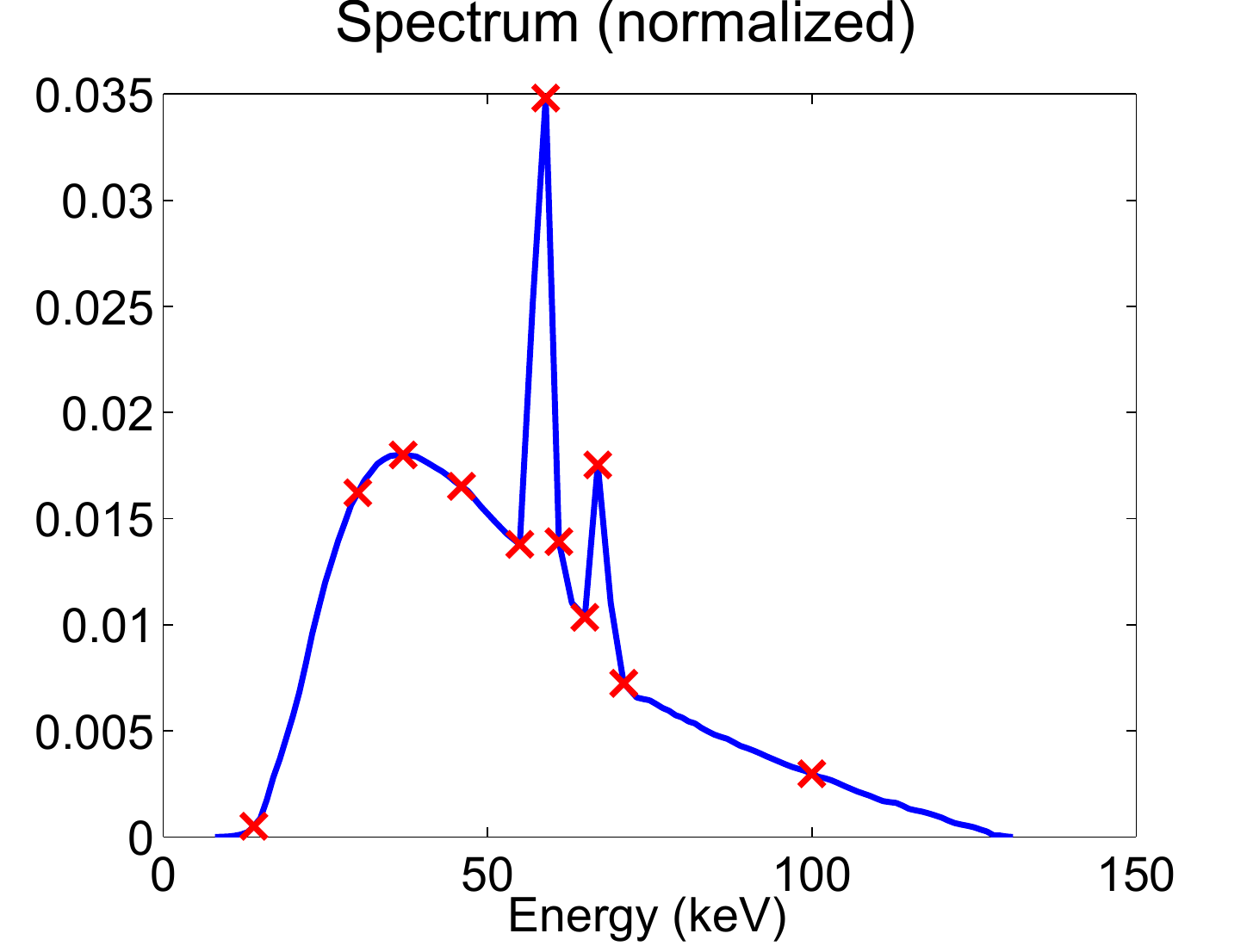} &\includegraphics[width=0.45 \linewidth]{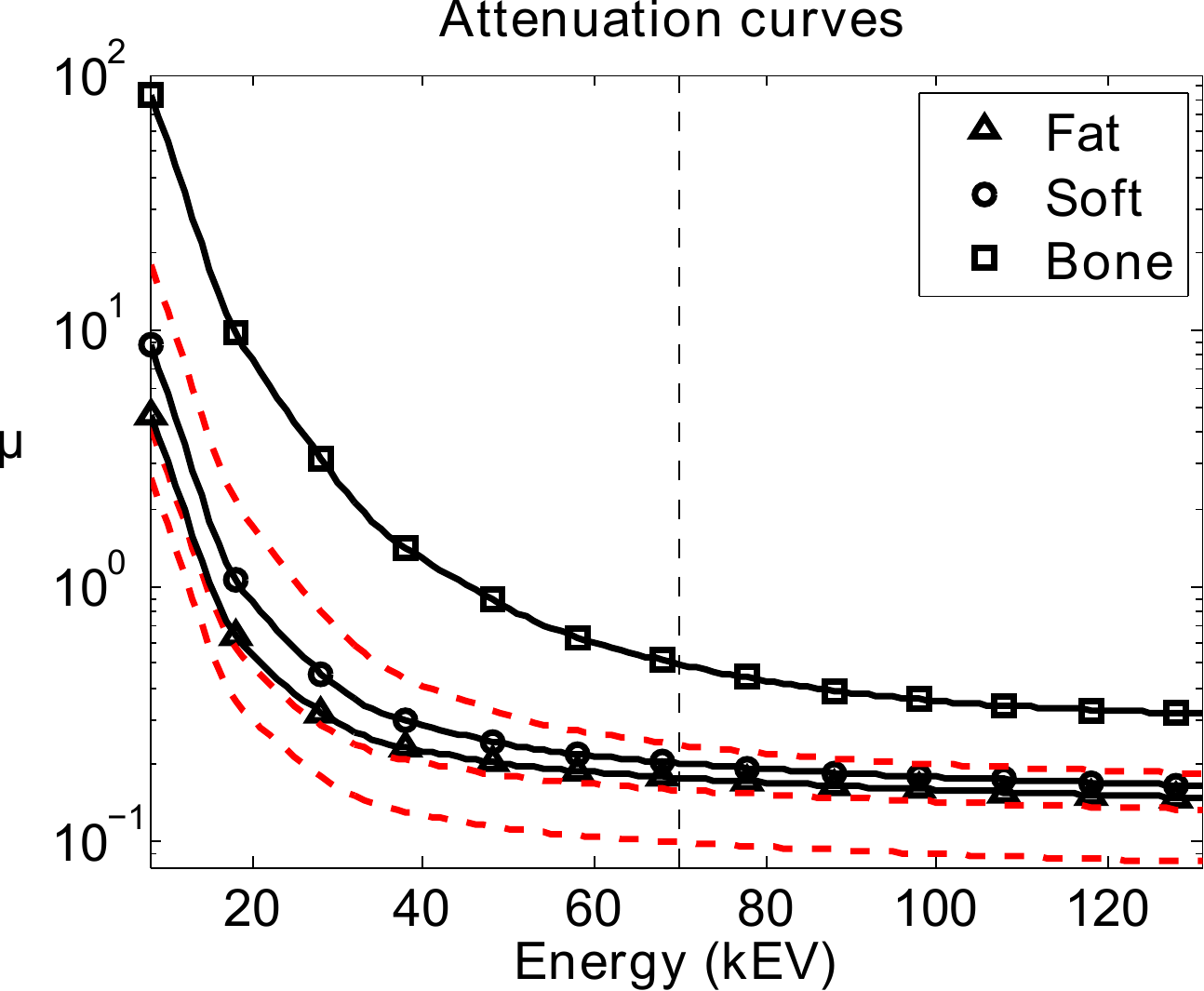}
\end{tabular}
\caption{X-ray spectrum and attenuation curves used for the polyenergetic experiment. Left: Continuous spectrum (blue line) and discrete energies (red crosses) used for the summation. The spectrum has been normalized to have an integral of 1. Right: Attenuation curves for the base materials as well as the interpolated curves for values of $t=0.1$, $0.16$ and $0.24$ (red dashed lines). The reference energy of 70 keV is indicated by the dashed black line.}\label{F:spectrum}
\end{figure}

\begin{figure}
\begin{tabular}{cc}
\includegraphics[width=0.45 \linewidth]{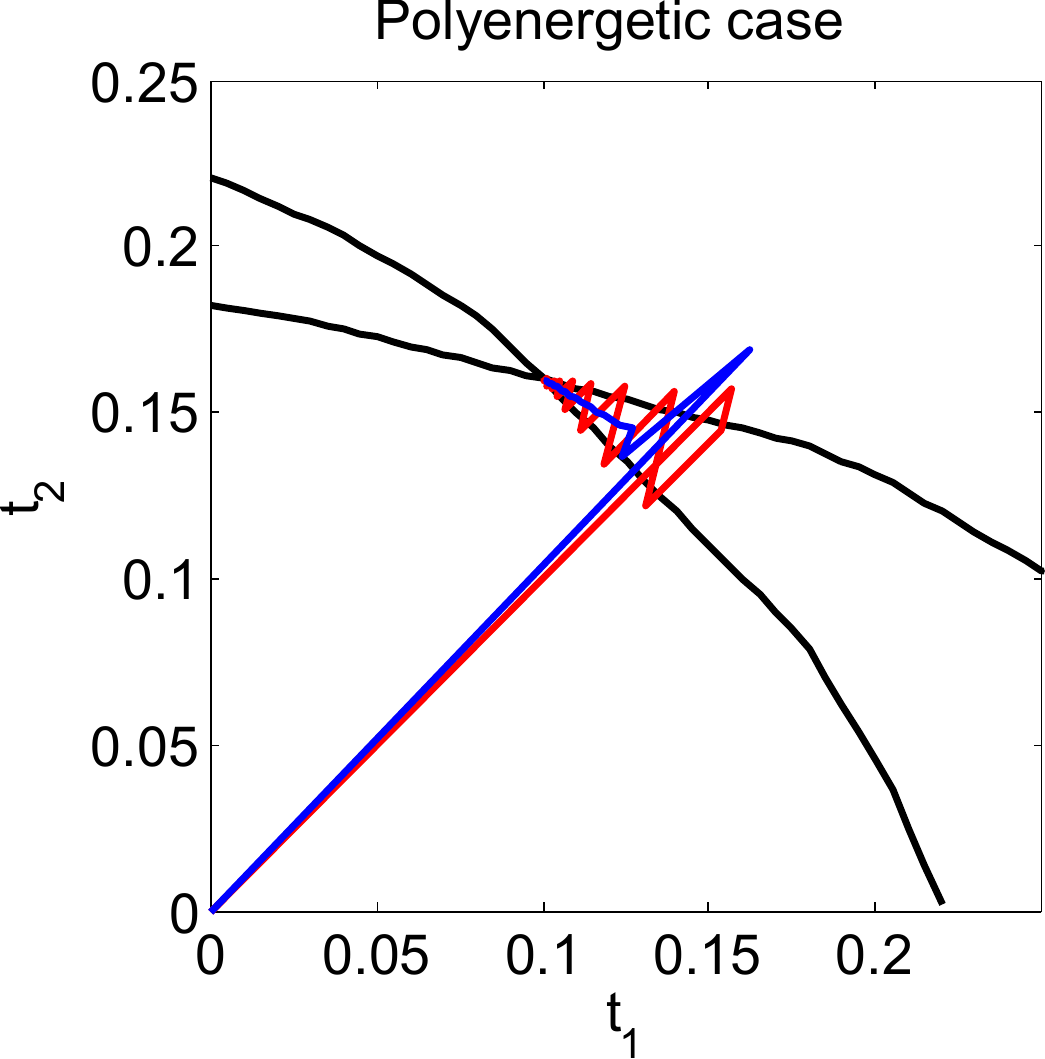} \includegraphics[width=0.45\linewidth]{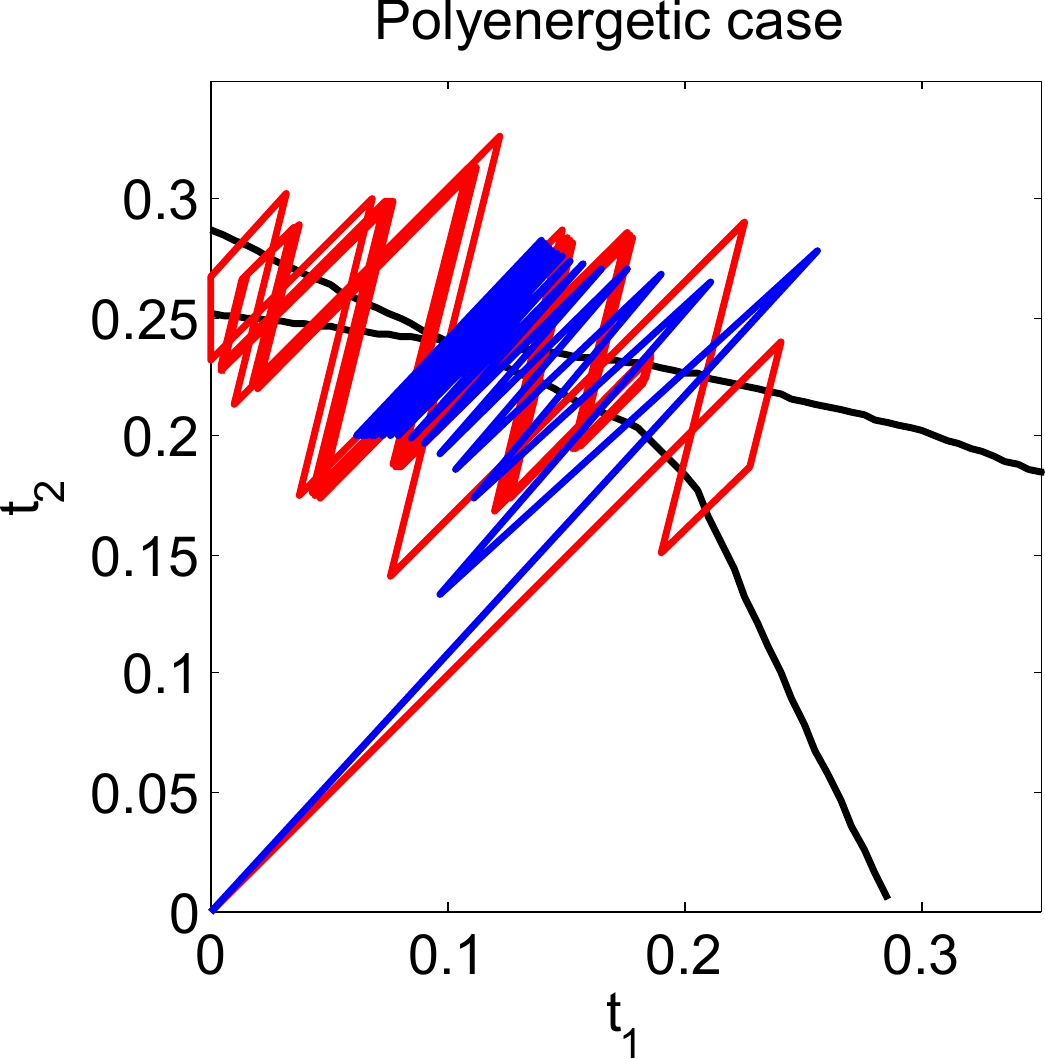}
\end{tabular}
\caption{Convergence of pART and pSART for two polyenergetic experiments. The black curves correspond to the nonlinear system of equations representing the polyenergetic mode, while the red and blue lines correspond to the sequence of iterates generated by the pART and pSART algorithms, respectively. In the first experiment (left plot) the LAC values at the reference energy are $t_1 = 0.1$ and $t_2 = 0.16$, while in the left plot, $t_2 = 0.24$.}\label{F:ART_poly1}
\end{figure}

We now consider the case of a polyenergetic x-ray spectrum. The LACs of the two pixels at the reference energy determine the LACs at all other energies $\varepsilon$ according to~(\ref{E:interp}). The discrete spectrum used in this experiment consists of eleven energy bins, obtained from a continuous 130 kVp spectrum generated using the Siemens Spektrum online tool~\cite{spektrum, BS97}. This spectrum and the attenuation curves for the two materials are shown in Fig.~\ref{F:spectrum}. The attenuation values for the reference materials were obtained from Ref.~\onlinecite{NIST}. Under the polyenergetic model, we obtain a nonlinear system of equations:
\begin{align*}
p_1 &= \sum_h I_h \exp \left( - \mu (t_1, \varepsilon_h) - \mu(t_2, \varepsilon_h) \right)  \\
	& \approx \exp(-0.314) \\
p_2 &= \sum_h I_h \exp \left( - 0.28 \mu (t_1, \varepsilon_h) - 1.13 \mu(t_2, \varepsilon_h) \right)  \\
	& \approx \exp(-0.253)
\end{align*}
The two beams undergo more attenuation than in the monoenergetic experiment (\ref{E:eqs1}) because the spectrum contains a higher proportion of X -rays with energies less than 70 keV. In Fig.~\ref{F:ART_poly1}, this nonlinear system of equations is illustrated with black curves. As either $t_1$ or $t_2$ increases, the LAC value of the pixel at lower energies increases rapidly, meaning that the coefficient in the other pixel must decrease rapidly to compensate. Thus the curves are bent. For the sake of illustration, we have implemented a polyenergetic version of ART (denoted pART) that is analogous to pSART. The red and blue lines show the progression of the pART and pSART iterations, respectively. For this experiment (left figure), one can see that both iterations converge to the solution of the nonlinear system of equations.

We now give a case where the iteration fails to converge. In the right figure of Fig.~\ref{F:ART_poly1}, the progress of the two iterations is shown for a slightly modified experiment, where $t_2$ was changed from 0.16 cm$^{-1}$ to 0.24 cm$^{-1}$. All other parameters of the experiment were the same as before. It is apparent that both iterations fail to converge to the solution; pART appears to exhibit chaotic behaviour about the solution, while pSART converges to a two-cycle. Even if the iteration is started very close to the solution $\tvec^* = (0.1, 0.24)^T$, both pART and pSART diverge.

\begin{figure}
\begin{tabular}{cc}
\includegraphics[width=0.45 \linewidth]{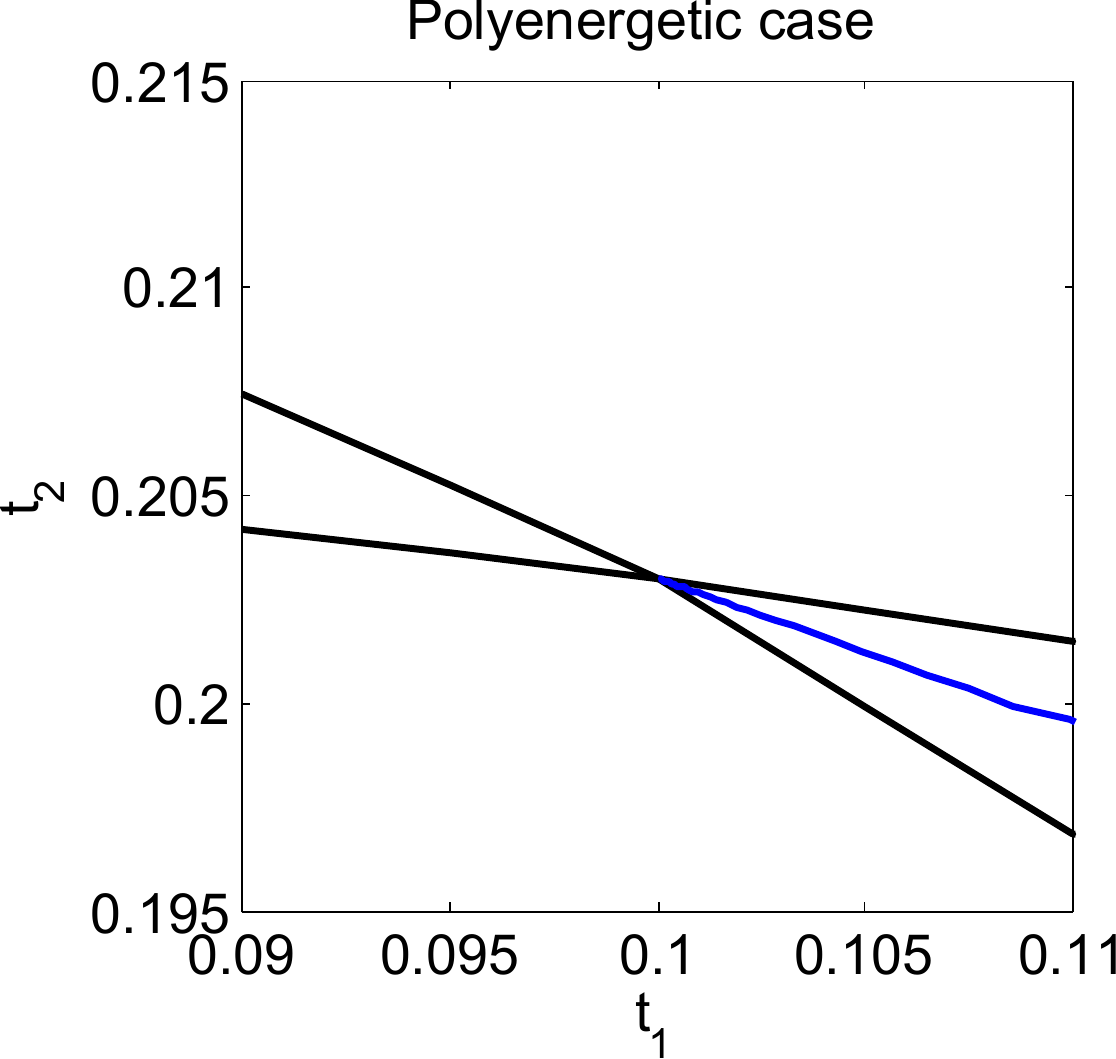} &\includegraphics[width=0.45 \linewidth]{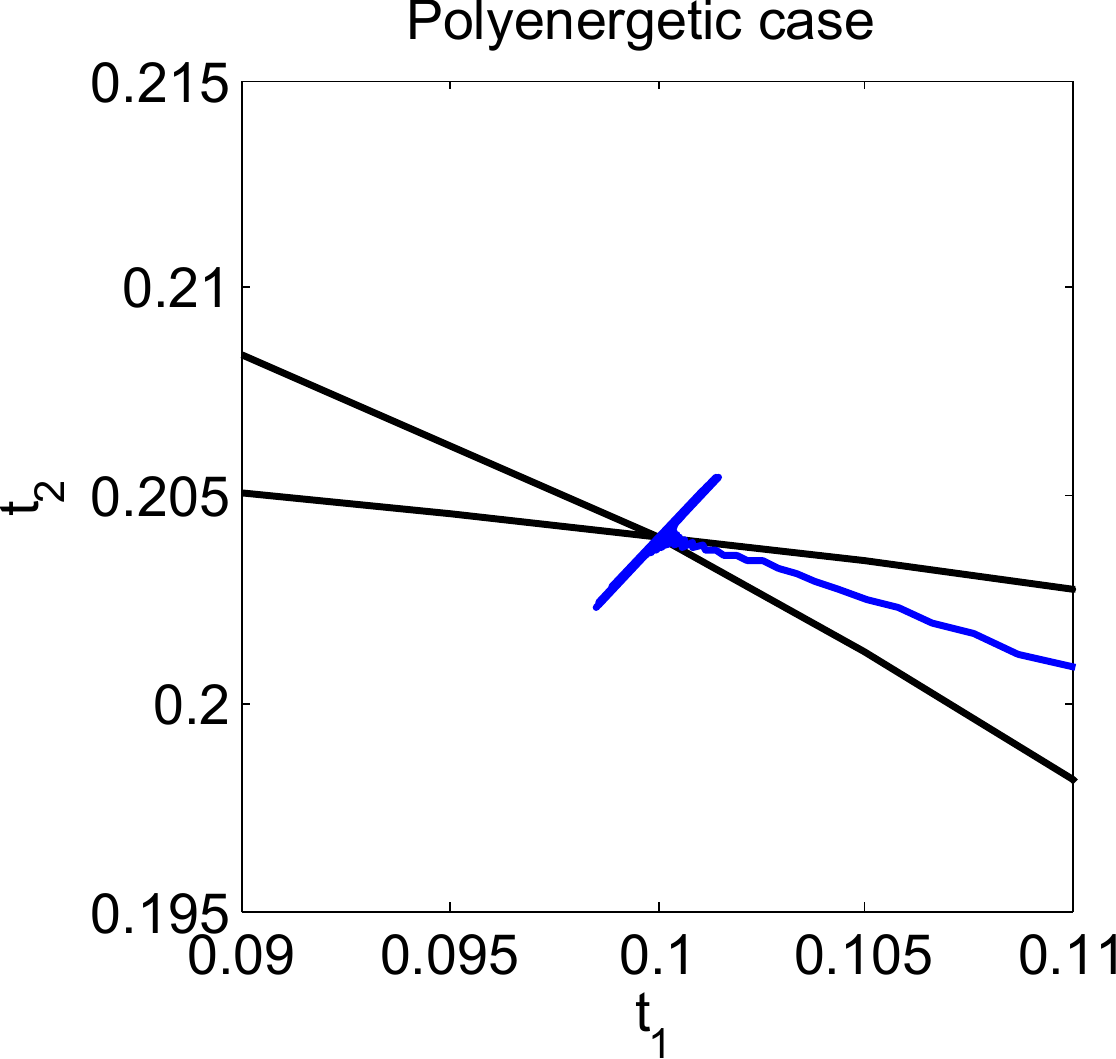}
\end{tabular}
\caption{pSART iterations for $t_2 = 0.203$ (left) and $t_2 = 0.204$ (right), narrowly centred on the solution. The first set of iterations converges to solution while the second converges to a two-cycle.}\label{F:ART_poly2}
\end{figure}

A direct computation of the $2 \times 2$ Jacobian matrix~(\ref{E:jac1}) for these two experiments reveals that the spectral radius of $J_F$, evaluated at the solution $\tvec^*$, is roughly 0.89 for the first case ($t_2 = 0.16$) and 1.02 for the second case ($t_2 = 0.24$). This explains why the first iteration converges, but not the second. Some further investigation reveals that this is due in large part to the discontinuities in $\displaystyle \frac{ \partial \mu}{\partial t}$. In our experiment the reference materials were air, fat, soft tissue and bone, with tabulated LAC values at the reference energy of 70 keV equal to 0, 0.1782, 0.2033 and 0.4948 cm$^{-1}$, respectively. Thus the partial derivative $\displaystyle \frac{ \partial \mu}{\partial t}$ has a larger value for $t_2 = 0.24$ (which lies between soft tissue and bone) than for $t_2 = 0.16$ (which lies between fat and air). In Fig.~\ref{F:ART_poly2} we show the result of two more experiments where $t_2$ was set to values of 0.203 and 0.204, which lie on either side of the tabulated value for soft tissue, where the derivative is discontinuous. The pSART iteration converges to the true solution in the first case but not in the second case, where it reaches a two-cycle between two points lying close to the true solution. Direct computation of the spectral radius confirms that it is equal to 0.87 in the first case and 1.16 in the second.

\begin{figure}
\includegraphics[width=0.95 \linewidth]{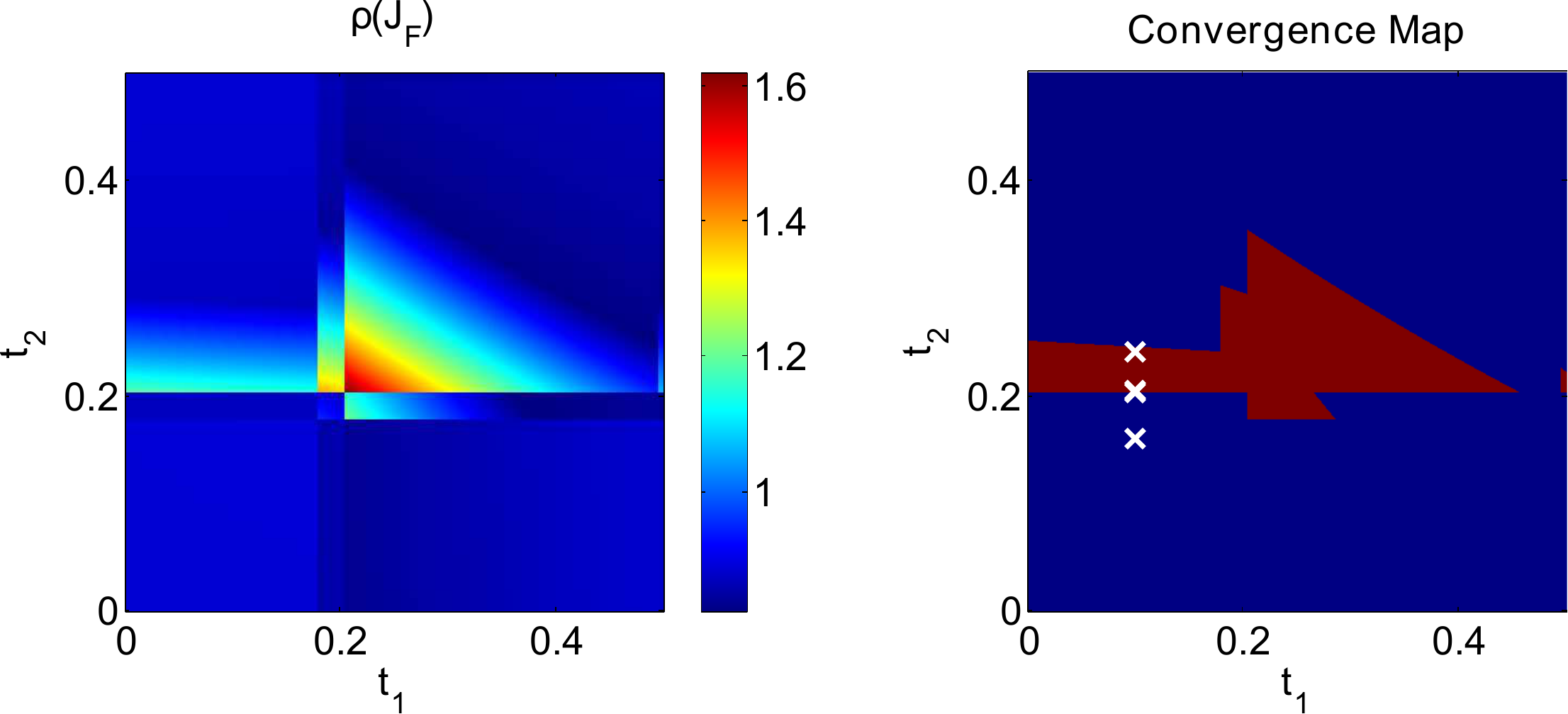} 
\caption{Spectral radius and convergence maps for the 2$\times$2 polyenergetic experiment. Left figure shows the spectral radius $\rho(J_F)$ as a function of $t_1$ and $t_2$; right figure is the equivalent binary map indicating whether the iteration converges (blue) or not (red). White crosses indicate the cases shown in Figs.~\ref{F:ART_poly1} and \ref{F:ART_poly2}.}\label{F:conv_map}
\end{figure}

Fig.~\ref{F:conv_map} gives a convergence map for this test case as a function of $t_1$ and $t_2$. The effect of the discontinuities in $\displaystyle \frac{ \partial \mu}{\partial t}$ is clearly visible in the discontinuities in $\rho(J_F)$ that occur at the values of 0.1782 -- the LAC value of fat at 70 keV -- and 0.2033, the LAC value of soft tissue at 70 keV. It is apparent that the iteration transitions between convergent and non-convergent states at values of $(t_1,t_2)$ that do not lie along these discontinuities as well. We note that the figure is not symmetric along the line $t_1=t_2$; for example, all of the test cases that have been considered in this section would converge if the values of $t_1$ and $t_2$ were interchanged. This figure is specific to the system matrix that arises from the ray paths illustrated in Fig.~\ref{F:exp1}, and would be different for other paths.


\section{Numerical experiments}
We have shown in the previous section that the pSART iteration is not guaranteed to converge in general, despite the success of the algorithm demonstrated in Ref.~\onlinecite{LS14b}. One possible explanation is that the system matrices that arise in CT imaging are typically quite sparse and structured, compared to the $2 \times 2$ matrix that was considered in the previous section. Thus, it is worth investigating whether the spectral radius of the Jacobian matrix of pSART is likely to exceed one for more realistic CT system matrices. 

\begin{figure}
\includegraphics*[width=0.6\linewidth]{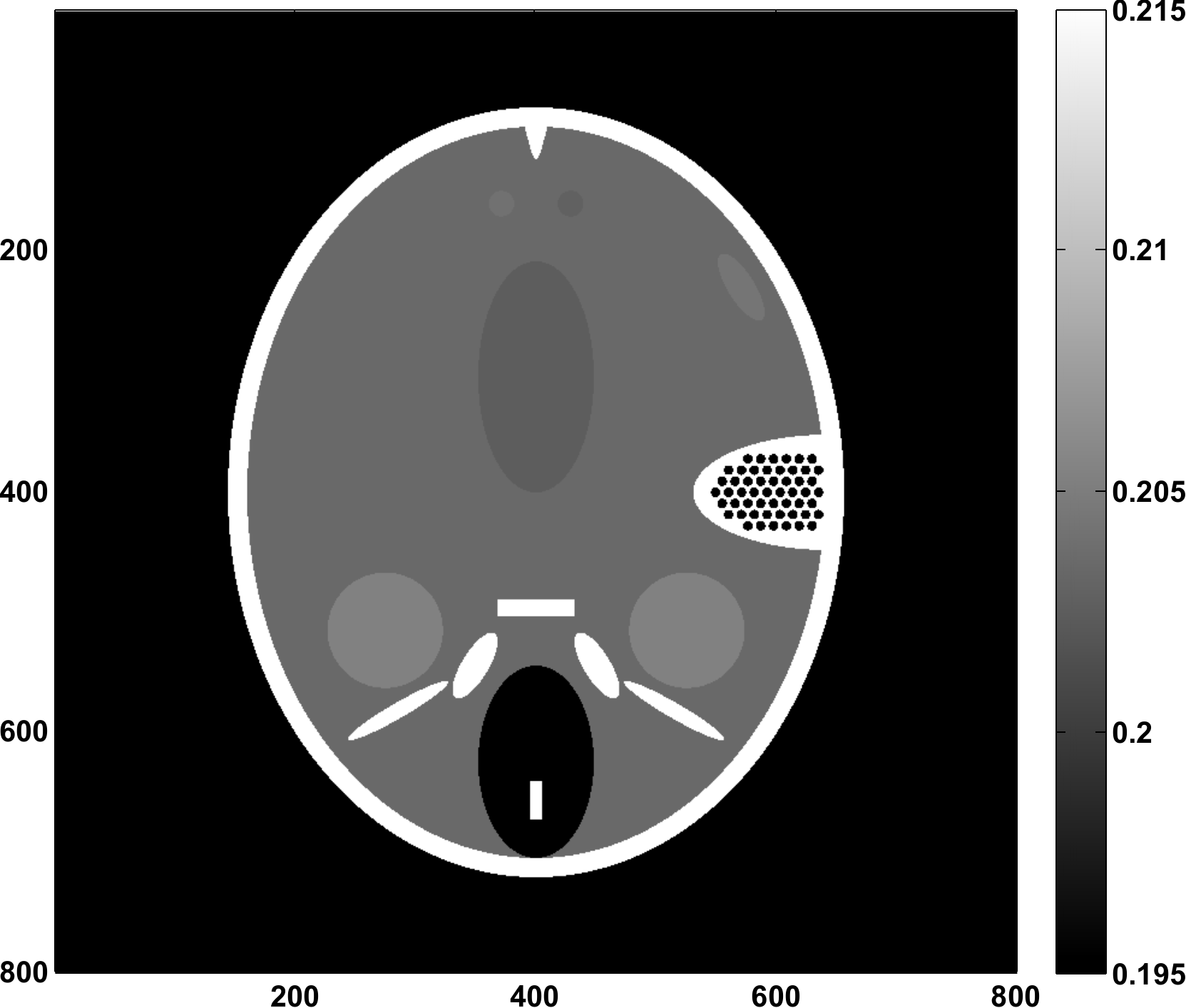}
\caption{FORBILD phantom for the case $N=800$. Color window is restricted to $[0.195, 0.215]$ to show low-contrast features; bony structures (white) have a LAC of 0.4948 at the reference energy of 70 keV.}\label{F:FORBILD}
\end{figure}

In the following numerical experiment we consider the problem of reconstructing an $N \times N$ pixel image for $N$ = 100, 200, 400 and 800. We simulate parallel beam data acquired at $m$ equally spaced views over 180$^\circ$, with $m$ = 180, 360, 720 and 1440, respectively. Forward projection (multiplication by $A$) is implemented using the {\tt radon} command in Matlab, while backward projection (multiplication by $A^T$) uses {\tt iradon} with no filtering. Since $J_F(\tvec^*)$ depends on the object $\tvec^*$ that we wish to reconstruct, we must consider a specific object to analyze the convergence of pSART. We use an $N \times N$ slice of the FORBILD numerical head phantom~\cite{YNDW12}, which consists of bone and soft tissue, and includes some low-contrast features. An image of the phantom for the case $N=800$ is shown in Fig.~\ref{F:FORBILD}. 

With these elements in place, we can approximate the spectral radius of the matrices associated with the SART iteration~(\ref{E:SART2}) and the Jacobian matrix of the pSART iteration (\ref{E:jac1}, \ref{E:jac2}). To approximate the spectral radius, we use the power iteration (see for instance Ref.~\onlinecite{TB97}), which provides iterative estimates of the largest eigenvalue of a matrix $M$ and the associated eigenvector, denoted by $\lambda^{(k)}$ and $\xvec^{(k)}$, respectively. The iterations were started from a random unit vector $\xvec^{(0)}$ and run until the the largest element of the residual, $\displaystyle \Vert M \xvec^{(k)} - \lambda^{(k)} \xvec^{(k)} \Vert_\infty$, was less than $\displaystyle 10^{-4}$.


Table~\ref{T:exp1} shows the computed approximations to $\rho(T)$ and $\rho(J_F(\tvec^*))$ for the different test cases, along with the number of power iterations required to obtain the estimate. It is apparent that there is virtually no difference between the spectral radius of the SART iteration matrix, $T$, and that of the the pSART Jacobian matrix, $J_F(\tvec^*)$, in any of the studied cases. In no instances did the computed estimate of $J_F(\tvec^*)$ ever exceed one, which would cause the iteration to diverge in the neighbourhood of the solution. We conclude that for more realistic CT imaging scenarios, the pSART iteration is, at the very least, likely to exhibit local convergence in the neighbourhood of the solution, with a comparable rate of convergence to SART. We note, however, that the spectral radius of both the SART and pSART iterations is very close to one, indicating that the convergence will be slow in the neighbourhood of the solution. The convergence can likely be accelerated by using subsets of the projection data, as was proposed in the original SART algorithm~\cite{AK84}.

\begin{table}
\caption{Results of the power iteration for the matrix $T$ used for the SART iteration, and the Jacobian matrix of the pSART iteration, for different phantom sizes. $N$ is the dimension the image to be reconstructed and $m$ is the total number of views. $\rho(T)$ is the spectral radius of the matrix $T$, and $\rho(J_F(\tvec^*))$ is the spectral radius of the Jacobian of pSART iteration, evaluated for the $N \times N$ FORBILD phantom. The number of iterations run for each power iteration are displayed below the estimate of the spectral radius.}\label{T:exp1}
\begin{tabular}{rrrr}
\hline
$N$	&$m$ 			&$\rho(T)$ 	 	&$\rho(J_F(\tvec^*))$	\\
\hline
100		&180		& 0.999957    	&0.999960  	\\
\smallskip		&		&2535 its		 &2802 its \\ 
200		&360		& 0.999957 		&0.999958 			 \\
\smallskip		&		&2899 its		&2880 its	\\	
400		&720		&0.999952 		&0.999954  		 \\
\smallskip		&		&3475 its		&3544 its	\\
800		&1440		&0.999933 	&0.999935			 \\
		&		&3544 its 	&3595 its\\
\hline
\end{tabular}
\end{table}

This analysis establishes only local convergence in the neighbourhood of the solution. Global convergence (i.e. from an arbitrary initial estimate) is more difficult to establish in general, and it is not obvious whether there exist conditions on the system matrix $A$, beam spectrum, choice of reference energy, etc. to guarantee global convergence of pSART. This is fairly typical of nonlinear fixed-point iterations, Newton's method being a well-known example. The method does appear to be fairly robust with respect to the choice of starting point, however, as the images in Ref.~\onlinecite{LS14b} were produced from an initial estimate consisting only of zeros; our own experience with the algorithm confirms that this choice of initial estimate works well.

Additionally, our analysis assumes the existence of an exact solution, whereas in practice the system is likely to be inconsistent due to factors such as noisy measurements and model mismatch. As with any reconstruction algorithm, model mismatch will produce artifacts in the reconstructed image; Ref.~\onlinecite{LS14b} includes some experiments quantifying the effect of spectrum error. In the presence of noisy data, our experience indicates that pSART exhibits the ``semi-convergence'' behaviour typical of other iterative methods (see e.g. Ref.~\onlinecite{N86}, p.89); namely, that the algorithm initially converges towards the solution, but the image eventually deterioriates with further iterations due to the effects of noise. As noted in Ref.~\onlinecite{LS14b}, this problem could potentially be addressed with the use of statistical modeling or edge-preserving regularization.

\section*{Conclusions}
In this paper we have analyzed the convergence of a recently proposed polyenergetic SART (pSART) algorithm. We show that the spectral radius of the Jacobian of the nonlinear pSART iteration may be larger than one in some cases. Thus the method is not mathematically guaranteed to converge to a solution of the nonlinear system of polyenergetic equations, in general. For system matrices of the type encountered in CT imaging, however, our empirical results indicate that the spectral radius of the Jacobian matrix, evaluated for a prototypical head phantom, is essentially the same as the spectral radius of the matrix corresponding to the convergent SART iteration. Thus in practice it seems that the method is likely to converge at roughly the same rate as SART.

\section*{Acknowledgments}

The author thanks Adel Faridani (Oregon State University) and Yuan Lin (Duke University) for helpful discussion on this paper, and  Federico Poloni (Universit\`{a} di Pisa) for helpful suggestions on the proof of Lemma A.2 in the Appendix.

\section*{Appendix}

We provide a short proof of convergence for SART subject to a condition on the rank of $A$. Convergence of SART has been established previously in Refs. \onlinecite{CE02} and \onlinecite{JW03}. This proof is somewhat shorter and is intended to complement the analysis we have presented for pSART.

Let $T = I-DA^TMA,$ where $A$ is the $m \times n$  system matrix and $D$ and $M$ are defined in~(\ref{E:diagmatrices}). We assume that $A$ has rank $n$, meaning that there are as many linearly independent equations as there are unknowns. This implies that $m \geq n$, i.e. that there are at least as many measurements as there are unknowns. The SART iteration then has the form
\begin{align*}
\xvec^{(k+1)} = T\xvec^{(k)} + \cvec
\end{align*}

We first prove two lemmas.

\smallskip

\noindent {\bf Lemma A.1:} {\em  Let $W = DA^TMA$. Then, $\rho(W) \leq 1$, with equality if all elements of $A$ are positive}.

\noindent {\em Proof:}  A direct calculation gives the  $(i,j)$th element of $W$ as
\begin{align*}
w_{ij} &= \frac{1}{\beta_i} \sum_{k=1}^m \frac{a_{ki} a_{kj}}{\gamma_k}.
\end{align*}
It follows that the sum of the absolute values in row $i$ of $W$ is:
\begin{align*}
\sum_{j=1}^n |w_{ij}| 								&= \sum_{j=1}^n \frac{1}{\beta_i} \biggl|\sum_{k=1}^m \frac{a_{ki} a_{kj}}{\gamma_k} \biggl| \\
												&\leq \frac{1}{\beta_i} \sum_{k=1}^m \frac{|a_{ki}|}{\gamma_k} \sum_{j=1}^n |a_{kj}| \\
												&=  \frac{1}{\beta_i} \sum_{k=1}^m \frac{|a_{ki}|}{\gamma_k} \gamma_k \\
												&= \frac{1}{\beta_i} \sum_{k=1}^m |a_{ki}| \\
												&= \frac{1}{\beta_i} \beta_i \\
												&= 1
\end{align*}

Since the max norm of a matrix, $\Vert \cdot \Vert_{\infty}$, is equal to the maximum row sum, and the spectral radius of a matrix cannot be greater than the max norm, it must be true that
\begin{equation*}
\rho(W) \leq \Vert W \Vert_{\infty} \leq 1.
\end{equation*}
When all elements of $A$ are positive (which is the case in tomographic applications), every row of $W$ sums exactly to 1, and so $\lambda = 1$ is an eigenvalue of $W$ (with the associated eigenvector consisting of all ones), and $\rho(W) = 1$. 

On its own, this lemma only tells us that the eigenvalues of $W$ have magnitude less than 1. We also need the following result:

\smallskip

\noindent {\bf Lemma A.2:} {\em All eigenvalues of $W$ are positive real numbers.}

\noindent {\em Proof:} This Lemma is an application of Theorem 7.6.3 from Ref.~\onlinecite{HJ85}. We first remind the reader of two definitions:
\begin{enumerate}
\item Two $n \times n$ matrices $A$ and $B$ are {\em similar} if there exists an invertible matrix $P$ such that $A = PBP^{-1}$. If $A$ and $B$ are similar, then they have exactly the same eigenvalues.
\item Two real $n \times n$ matrices $A$ and $B$ are {\em congruent} if there exists an invertible matrix $P$ such that $A = PBP^T$. If $A$ and $B$ are symmetric, then Sylvester's law of inertia states that they have the same inertia, meaning the same number of positive, negative, and zero eigenvalues. (Recall that the eigenvalues of a symmetric matrix are always real numbers).
\end{enumerate}

Now, let $V = A^TMA$. Then, V is a symmetric positive semidefinite $n \times n$ matrix, since it can be written as the product of a matrix and its transpose: $\displaystyle V = (M^{\frac{1}{2}}A)^T (M^{\frac{1}{2}}A)$. Furthermore, since $A$ has rank $n$ and $\displaystyle M^{\frac{1}{2}}$ is a diagonal matrix, $\displaystyle M^{\frac{1}{2}}A$ has rank $n$ and so does $V$. It follows that $V$ is invertible, and so zero is not an eigenvalue of $V$. Thus all eigenvalues of $V$ are positive. We then have the following two results:
\begin{enumerate}
\item $W = DV$ is similar to $\displaystyle D^{\frac{1}{2}} V D^{\frac{1}{2}}$, since $\displaystyle D^{\frac{1}{2}} V D^{\frac{1}{2}} = D^{-\frac{1}{2}} D V D^{\frac{1}{2}}$.
\item $D^{\frac{1}{2}} V D^{\frac{1}{2}}$ is congruent to $V$.
\end{enumerate}
The first result implies that the eigenvalues of $W$ are the same as the eigenvalues of $\displaystyle D^{\frac{1}{2}} V D^{\frac{1}{2}}$, while the second implies that the eigenvalues of this matrix must have the same signs as the eigenvalues of $V$. So, since all eigenvalues of $V$ are positive, all eigenvalues of $W$ must be positive as well. 
\smallskip

\noindent {\bf Theorem A.3:} The matrix $T$ satisfies $\rho(T) < 1$, and hence the SART iteration converges.

\noindent {\em Proof:} Lemmas A.1 and A.2 prove that any eigenvalues $\lambda$ of $W$ satisfy $0 < \lambda \leq 1$. Thus the eigenvalues of $T = I-W$ satisfy $0 \leq \lambda < 1$, and so $\rho(T) < 1$.

\end{document}